\documentclass{article}

\usepackage{graphicx}
\usepackage{amsmath}
\usepackage{amsfonts}
\usepackage{amssymb}

\newtheorem{theorem}{Theorem}

\newtheorem{corollary}{Corollary}

\newtheorem{lemma}{Lemma}

\newtheorem{proposition}{Proposition}

\newenvironment{proof}[1][Proof]{\textbf{#1.} }{\ \rule{0.5em}{0.5em}}

\begin{document}
\title{Non-semisimple Lie algebras with Levi factor $\frak{so}(3), \frak{sl}(2,\mathbb{R})$ and their invariants}

\author{Rutwig Campoamor\\Departamento de Geometr\'{\i}a y Topolog\'{\i}a\\Fac. CC. Matem\'{a}ticas\\Universidad Complutense\\28040 Madrid ( Spain )\\rutwig@nfssrv.mat.ucm.es}
\date{}
\maketitle

\begin{abstract}
We analyze the number $\mathcal{N}$ of functionally independent generalized
Casimir invariants for non-semisimple Lie algebras $\frak{s}\overrightarrow{%
\oplus }_{R}\frak{r}$ with Levi factors isomorphic to $\frak{so}\left(
3\right) $ and $\frak{sl}\left( 2,\mathbb{R}\right) $ in dependence of the
pair $\left( R,\frak{r}\right) $ formed by a representation $R$ of $\frak{s}$
and a solvable Lie algebra $\frak{r}$. We show that for any dimension $n\geq
6$ there exist Lie algebras $\frak{s}\overrightarrow{\oplus }_{R}\frak{r}$
with non-trivial Levi decomposition such that $\mathcal{N}\left( \frak{s}%
\overrightarrow{\oplus }_{R}\frak{r}\right) =0$.   
\end{abstract}

\section{Introduction}

The important role played by invariant theory in Physics has been
recognized long ago. Electroweak interactions and quantum numbers in the study
of particle states are based on the concept of symmetry, and their invariants
provide fundamental information. Among the various types of symmetry, dynamical
ones constitute one of the more important cases, as shown by Gell-Mann and
Ne'eman in their hadron classification \cite{Ge}. The analysis of the group
$SU\left(  3\right)  $ resulted in the prediction of new particles whose mass
could be derived from the invariants of the group. The invariants of Lie
algebras have also shown their efectiveness in the description of Hamiltonians
\cite{Fo}, the labelling
of irreducible representations or the study of coadjoint orbits \cite{Ro,Ki}.
Other important applications of invariants arise in their combination with the
theories of Lie algebra contractions, deformations and rigidity
\cite{Cam,AC,Car,Cam2}. For example, all kinematical algebras are related by a
contraction procedure, which has allowed a further analysis of these algebras
\cite{Car,Lo,Ba}. The interest of invariants of rigid Lie algebras is fully
justified by the fact that semisimple Lie algebras are rigid. The invariants of semisimple Lie algebras constitute a classical problem, and it is the only case which has been solved in a satisfactory manner. The invariants of solvable Lie algebras are only studied for specific classes, as they do not underly to a structure theory like the classical algebras. What refers to the  Lie algebras with non-trivial Levi decomposition, invariants are known for physically important algebras, like the special affine algebras $\frak{sa}(n,\mathbb{R})$, the kinematical Lie algebras and their subalgebras.

\medskip

A formula for the number $\mathcal{N}\left(  \frak{g}\right)  $ of
functionally independent invariants of the coadjoint representation of a Lie
algebra $\frak{g}$ was given by Beltrametti and Blasi \cite{Be} and Pauri and
Prosperi \cite{Pa} in the mid 1960s. This fact reduces the computation of this number to the
determination of the rank of a skew-symmetric matrix $A\left(  \frak{g}\right)  $ whose
entries correspond to the Lie brackets of $\frak{g}$. With some effort, this formula can be used to show that the number of invariants of semisimple Lie algebras coincides with
its rank \cite{Ra}. Moreover, it proves that for direct sums $\frak{g}_{1}\oplus
\frak{g}_{2}$ of Lie algebras the number $\mathcal{N}\left(  \frak{g}%
_{1}\oplus\frak{g}_{2}\right)  $ is $\mathcal{N}\left(  \frak{g}_{1}\right)
+\mathcal{N}\left(  \frak{g}_{2}\right)  $. One can ask whether for
semi-direct sums $\frak{g}=\frak{s}\overrightarrow{\oplus}_{R}\frak{r}$,
$\frak{s\,}$\ being the Levi factor of $\frak{g}$, $R$ a representation of
$\frak{s}$ and $\frak{r}$ the maximal solvable ideal (called radical) of $\frak{g}$ some
formula exists which allows to express $\mathcal{N}\left(  \frak{s}%
\overrightarrow{\oplus}_{R}\frak{r}\right)  $ in terms of $\mathcal{N}\left(
s\right)  ,\mathcal{N}\left(  \frak{r}\right)  $ and some quantity related to
the representation $R$. The motivation of this problem lies in the study of
the special affine algebras $\frak{sa}\left(  n,\mathbb{R}\right)  $, which
are a semidirect sum of the simple Lie algebra $\frak{sl}\left(
n,\mathbb{R}\right)  $ and an $n$-dimensional abelian Lie algebra \cite{PeN}. These
algebras are known to have only one invariant (which turns out to be a Casimir
operator), which shows that the representation plays a crucial role in the
semidirect product, and that in principle the existence of a formula
expressing the number of invariants in terms of the factors does not exist.
The main reason for its nonexistence lies in the distinct possibilities of
choice for radicals $\frak{r}$ \ for a fixed representation of $\frak{s}$. The
question that arises naturally in this context is if there exist Lie algebras
$\frak{s}\overrightarrow{\oplus}_{R}\frak{r}$ \ with non-trivial Levi
decomposition (i.e. $\frak{s}\neq0$ and $\left[  \frak{s},\frak{r}\right]
\neq0$) such that $\mathcal{N}\left(  \frak{s}\overrightarrow{\oplus}%
_{R}\frak{r}\right)  =0$. \ 

In this work we show that such algebras exist for any dimension
$n\geq6$. Moreover, by considering the simple algebras $\frak{so}\left(
3\right)  $ and $\frak{sl}\left(  2,\mathbb{R}\right)  $, we analyze the number
$\mathcal{N}\left(  \frak{s}\overrightarrow{\oplus}_{R}\frak{r}\right)  $ for
various kinds of representations $R$ and solvable Lie algebras $\frak{r}$.

\medskip

Any Lie algebra $\frak{g}$ considered in this work is defined over the field $\mathbb{R}$ of real numbers. We convene that nonwritten brackets are either zero or obtained by antisymmetry. We also use the Einstein summation convention. Abelian Lie algebras of dimension $n$ will be denoted by $nL_{1}$.  

\section{Invariants of Lie algebras. The Beltrametti-Blasi formula}

The method to determine the invariants of a Lie algebra in terms of systems of
partial differential equations (PDEs) has become standard in the physical
literature \cite{AA,Pe}, and it is the one we will use here. Let $\left\{
X_{1},..,X_{n}\right\}  $ be a basis of $\frak{g}$ and $\left\{  C_{ij}%
^{k}\right\}  $ be the structure constants over this basis. We consider the
representation of $\frak{g}$ in the space $C^{\infty}\left(  \frak{g}^{\ast
}\right)  $ given by:
\begin{equation}
\widehat{X}_{i}=-C_{ij}^{k}x_{k}\partial_{x_{j}},
\end{equation}
where $\left[  X_{i},X_{j}\right]  =C_{ij}^{k}X_{k}$ $\left(  1\leq i<j\leq
n,\;1\leq k\leq n\right)  $. This representation is isomorphic to $ad\left(
\frak{g}\right)  $, and therefore satisfies the brackets $\left[  \widehat
{X}_{i},\widehat{X}_{j}\right]  =C_{ij}^{k}\widehat{X}_{k}$. The invariants
$F\left(  X_{1},..,X_{n}\right)  $ of $\frak{g}$:
\begin{equation}
\left[  X_{i},F\left(  X_{1},..,X_{n}\right)  \right],
\end{equation}
are found by solving the system of linear first order partial differential
equations:%
\begin{equation}
\widehat{X}_{i}F\left(  x_{1},..,x_{n}\right)  =-C_{ij}^{k}x_{k}%
\partial_{x_{j}}F\left(  x_{1},..,x_{n}\right)  =0,\;1\leq i\leq n.
\end{equation}
and then replacing the variables $x_{i}$ by the corresponding generator
$X_{i}$ (possibly after symmetrizing). In recent years new algorithms to solve
system $\left(  {3}\right)  $ have been developed, which simplify the
calculation in some cases \cite{Pe}. A maximal set of functionally independent
solutions of $\left(  {3}\right)  $ will be called a fundamental set of
invariants. Polynomial solutions of system $\left(  {3}\right)  $ are therefore
polynomials in the generators which commute with $\frak{g}$, thus correspond
to the well known Casimir operators \cite{AA}. The system does not impose
additional conditions which imply that the solutions are polynomials, so that
a non-polynomial solution will be called, in some analogy with the classical
case, a generalized Casimir invariant or simply an invariant of $\frak{g}$. If
$F$ reduces to a constant we say that the invariant is trivial. In the case of
semisimple Lie algebras, the solutions found are in fact Casimir operators,
and the number of functionally independent invariants is given by the
dimension of its Cartan subalgebra. However, for non-semisimple Lie algebras
there is no reason to suppose that only the polynomial invariants are of
physical interest. A classical example for a Hamiltonian being a nonlinear
function of the Casimir operators was described by Pauli in \cite{Pau}.

\medskip

Another important task is to find the maximal number $\mathcal{N}\left(
\frak{g}\right)  $ of functionally independent solutions of $\left(
{3}\right)  $. For the case of the classical groups this number depends only on
the dimension of a Cartan subalgebra, while for solvable Lie algebras no such
general formula exists \cite{Cam}. However, for a fixed algebra, this number
can be described in terms of the dimension and a certain matrix associated to
the commutator table. More specifically, denote by $A\left(  \frak{g}\right)
$ the matrix whose $\left(  i,j\right)  $-entry is the bracket $\left[
X_{i},X_{j}\right]  $. Such a matrix has necessarily even rank.\ Then
$\mathcal{N}\left(  \frak{g}\right)  $ is given by
\begin{equation}
\mathcal{N}\left(  \frak{g}\right)  =\dim\,\frak{g}-\left\{rank\,A\left(
\frak{g}^{\prime}\right)\quad |\quad \frak{g}^{\prime}\simeq \frak{g}\right\}   .
\end{equation}
This formula was first described by Beltrametti and Blasi \cite{Be} and Pauri
and Prosperi \cite{Pa}. The number of polynomial solutions is generally lower
than $\mathcal{N}(\frak{g})$ , up to certain special classes of Lie algebras
(like semisimple and nilpotent) \cite{AA}.\newline 

Invariants of Lie algebras have been determined for some classes of
non-semisimple Lie algebras, like solvable Lie algebras in low dimensions
\cite{Nd,PWZ}, the kinematical Lie algebras \cite{Car} or the special affine Lie
algebras \cite{Pe}.

\medskip

We give an example to illustrate the general method of obtaining the invariants.

\medskip

Let $\frak{s}=\frak{so}\left(  3\right)  $ and consider the representation
$R=ad\,\ \frak{so}\left(  3\right)  $. Let us suppose that the radical of the six
dimensional Lie algebra $\frak{s}\overrightarrow{\oplus}_{R}\frak{r}$ \ is
the three dimensional abelian algebra $3L_{1}$. The algebra $\frak{s}\overrightarrow{\oplus}%
_{R}\frak{r}$ \ is of interest for multidimensional extensions of the Bianchi
type-IX cosmology \cite{Tu}, and the corresponding vacuum Einstein field
equations have been solved in \cite{Tu}. Indeed this is the simplest
embedding of a Bianchi type-IX algebra in an algebra with non-trivial Levi
decomposition \cite{De}. It can easily be verified that $\frak{s}\overrightarrow{\oplus}_{R}\frak{r}$ 
satisfies $\mathcal{N}\left(  \frak{s}\overrightarrow{\oplus}%
_{R}\frak{r}\right)  =2$. The invariants are solution of the system:%

\begin{equation}
\left.
\begin{array}
[c]{l}%
\widehat{X}_{1}F=\left(  -x_{3}\partial_{x_{2}}+x_{2}\partial_{x_{3}}%
-x_{6}\partial_{x_{5}}+x_{5}\partial_{x_{6}}\right)  F=0\\
\widehat{X}_{2}F=\left(  x_{3}\partial_{x_{1}}-x_{1}\partial_{x_{3}}%
+x_{6}\partial_{x_{4}}-x_{4}\partial_{x_{6}}\right)  F=0\\
\widehat{X}_{3}F=\left(  -x_{2}\partial_{x_{1}}+x_{1}\partial_{x_{2}}%
-x_{5}\partial_{x_{4}}+x_{4}\partial_{x_{5}}\right)  F=0\\
\widehat{X}_{4}F=\left(  -x_{6}\partial_{x_{2}}+x_{5}\partial_{x_{3}}\right)
F=0\\
\widehat{X}_{5}F=\left(  x_{6}\partial_{x_{1}}-x_{4}\partial_{x_{3}}\right)
F=0\\
\widehat{X}_{6}F=\left(  -x_{5}\partial_{x_{1}}+x_{4}\partial_{x_{2}}\right)
F=0
\end{array}
\right\}  ,
\end{equation}
Since the equations $\left\{  \widehat{X}_{i}F=0\right\}  _{i=4,5,6}$ do not
depend on $\partial_{x_{i}}F$ for $i=4,5,6$, we can extract the following
system from $\left(  {5}\right)  :$%
\begin{equation}
\left.
\begin{array}
[c]{l}%
\widehat{X}_{1}^{\prime}F=\left(  -x_{6}\partial_{x_{5}}+x_{5}\partial_{x_{6}%
}\right)  F=0\\
\widehat{X}_{2}^{\prime}F=\left(  x_{6}\partial_{x_{4}}-x_{4}\partial_{x_{6}%
}\right)  F=0\\
\widehat{X}_{3}^{\prime}F=\left(  -x_{5}\partial_{x_{4}}+x_{4}\partial_{x_{5}%
}\right)  F=0
\end{array}
\right\}  ,
\end{equation}
which has the solution $I_{1}=x_{4}^{2}+x_{5}^{2}+x_{6}^{2}$. Now, as the rank
of the coefficients matrix corresponding to this subsystem is two, the other
solution of $\left(  {5}\right)  $ will depend also on $x_{1},x_{2},x_{3}$.
This invariant can be chosen as $I_{2}=x_{1}x_{4}+x_{2}x_{5}+x_{3}x_{6}$. The
important fact about this example is the solution found extracted from the
subsystem $\left(  {6}\right)  $. In a following section we will see that this
is not casual, but a property that holds in general.

\section{Semidirect sums of Lie algebras}

The classification of Lie algebras is simplified in some manner by the Levi
decomposition theorem, which states that any Lie algebra is essentially formed
from a semisimple Lie algebra $\frak{s}$ called the Levi factor of $\frak{g}$
and a maximal solvable ideal $\frak{r}$, called the radical \cite{RWZ}. Since the
latter is an ideal, the Levi factor $\frak{s}$ acts on $\frak{r}$, and there
are two possibilities for this action:%
\begin{equation*}
\begin{array}
[c]{l}%
\left[  \frak{s},\frak{r}\right]    =0\\
\left[  \frak{s},\frak{r}\right]    \neq0
\end{array}
\end{equation*}
If the first holds, then $\frak{g}$ is a direct sum $\frak{s}\oplus\frak{r}$,
whereas the second possibility implies the existence of a representation $R$ of
$\frak{s}$ which describes the action, i.e.,
\begin{equation}
\left[  x,y\right]  =R\left(  x\right)  .y,\;\forall x\in\frak{s},y\in\frak{r}%
\end{equation}
Unless there is no ambiguity, it is more convenient to write $\overrightarrow
{\oplus}_{R}$ instead of $\overrightarrow{\oplus}$, which is the common symbol
for denoting semidirect products. Since $\left(  {7}\right)  $ implies that the
radical is a module over $\frak{s}$, we have to expect severe restrictions on
the structure of the radical, while for direct sums any solvable Lie algebra
is suitable as radical.

\begin{proposition}
Let $\frak{s}$ be a semisimple Lie algebra and $R$ an irreducible
representation. If $\frak{s}\overrightarrow{\oplus}_{R}\frak{r}$ is the Levi
decomposition of a Lie algebra, then $\frak{r}$ is an abelian algebra.
\end{proposition}

The proof is immediate, since the Jacobi condition implies that the
ideals $\frak{r}^{\left(  0\right)  }:=\frak{r}$, $\frak{r}^{\left(  i\right)
}:=\left[  \frak{r}^{\left(  i-1\right)  },\frak{r}^{\left(  i-1\right)
}\right]  $ for $i\geq1$ are invariant by the action of $\frak{s}$. If $R$ is
irreducible, then either $\frak{r}^{\left(  1\right)  }=0$ or $\frak{r}%
^{\left(  1\right)  }=\frak{r}$, and since $\frak{r}$ is solvable,
$\frak{r}^{\left(  1\right)  }\neq\frak{r}$. Reasoning similarly, we can
easily deduce that the radical $\frak{r}$ is mapped into its maximal nilpotent
ideal $\frak{n}$ (usually called the nilradical of $\frak{r}$), from which the
following property follows:

\begin{proposition}
Let $\frak{s}\overrightarrow{\oplus}_{R}\frak{r}$ be a Levi decomposition. If
the representation $R$ does not posses a copy of the trivial representation,
then the radical $\frak{r}$ is a nilpotent Lie algebra.
\end{proposition}

This result is in some manner surprising, since it implies the existence of a copy of
the trivial representation whenever the radical is not nilpotent. Of course it does not imply that
a nilpotent Lie algebra cannot be the radical when the representation contains copies of the trivial 
representation.

\medskip

The Lie algebras having non-trivial Levi decomposition have been completely
classified up to dimension 8 \cite{Tu}. For dimensions 9 and 10 some partial results do also exist, mainly Levi factors isomorphic to rank one simple Lie algebras. Since the algebra $\frak{so}\left(  3\right)  $ is a real form
of $\frak{sl}\left(  2,\mathbb{C}\right)  $, the number of (real)
representations of $\frak{so}\left(  3\right)  $ is lower than for for
$\frak{sl}\left(  2,\mathbb{R}\right)  $ \cite{Iwa}, which implies the
existence of much more Lie algebras having the latter as Levi factor.

\begin{lemma}
Let $\frak{g}=\frak{s\oplus r}$. Then $\mathcal{N}\left(  \frak{g}\right)
=\mathcal{N}\left(  \frak{s}\right)  +\mathcal{N}\left(  \frak{r}\right)  $.
\end{lemma}

This is an obvious consequence of the Beltrametti-Blasi formula. Since the sum
is direct, we have that $\left[  \frak{s},\frak{r}\right]  =0$ and therefore
the rank of the matrix $A\left(  \frak{g}\right)  $ is the sum of the ranks of
$A\left(  \frak{s}\right)  $ and $A\left(  \frak{r}\right)  $. Now one can ask
what happens whenever we have a non-trivial Levi decomposition. Here no
apparent relation between the number of invariants of the Levi factor and the
radical and the number of invariants of the semidirect sum seems to exist. If
we consider the simple algebra $\frak{sl}\left(  2,\mathbb{R}\right)
=\left\{  X_{1},X_{2},X_{3}\;|\;\left[  X_{1},X_{2}\right]  =2X_{2},\;\left[
X_{1},X_{3}\right]  =-2X_{3},\;\left[  X_{2},X_{3}\right]  =X_{1}\right\}  $
and the representation $R=D_{\frac{1}{2}}\oplus D_{0}$, $D_{\frac{1}{2}}$
being the irreducible representation of highest weight $\lambda=1$, there are
two choices of $\frak{r}$ such that $\frak{sl}\left(  2,\mathbb{R}\right)
\overrightarrow{\oplus}\frak{r}$ is a six dimensional Lie algebra with
non-trivial Levi decomposition: either the 3-dimensional Heisenberg Lie
algebra $\frak{h}_{1}=\left\{  X_{4},X_{5},X_{6}\;|\;\left[  X_{4}%
,X_{5}\right]  =X_{6}\right\}  $ or the algebra $A_{3,3}=\left\{  X_{4}%
,X_{5},X_{6}\;|\;\left[  X_{i},X_{6}\right]  =X_{i},\;i=4,5\right\}  $ (see \cite{PWZ} for this notation). It is
a straightforward verification that $\mathcal{N}\left(  \frak{h}_{1}\right)
=\mathcal{N}\left(  A_{3,3}\right)  =1$, thus the formula, if existing, should
give the same value for both cases. Now the Lie algebra $\frak{sl}\left(
2,\mathbb{R}\right)  \overrightarrow{\oplus}_{R}\frak{h}_{1}$ admits two
(polynomial) invariants $I_{1}=x_{6}$ and $I_{2}=2x_{1}x_{4}x_{5}+4x_{2}%
x_{3}x_{6}+2x_{2}x_{5}^{2}-2x_{3}x_{4}^{2}+x_{1}^{2}x_{6}$, while the algebra
$\frak{sl}\left(  2,\mathbb{R}\right)  $ $\overrightarrow{\oplus}A_{3,3}$ has
no invariant. The conclusion is that the number of invariants will in general
be not expressible in terms of its factors. This example points out another
interesting fact: the existence of pairs $\left(  R,\frak{r}\right)  $ formed
by representations $R$ of a semisimple Lie algebra $\frak{s}$ and a solvable
Lie algebra $\frak{r}$ with structure of $\frak{s}$-module such that
\begin{equation}
\mathcal{N}\left(  r\right)  >0\quad {\rm and }\quad \mathcal{N}\left(  \frak{s}%
\overrightarrow{\oplus}_{R}\frak{r}\right)  =0.
\end{equation}
This also shows that it is not sufficient to determine the invariants of
solvable Lie algebras to have an overview of invariants of Lie algebras, implicitly
assumed in some early works. Thus the Levi decomposition theorem does not simplify
the determination of Casimir operators of Lie algebras, up to the case where
we obtain a direct sum. The next step is naturally to try the classification
of pairs $\left(  R,\frak{r}\right)  $ for fixed Levi factor $\frak{s}$ such
that $\left(  {8}\right)  $ holds. This problem cannot be solved since the
classification of solvable Lie algebras is probably not possible for
dimensions $n\geq7$ (the classification of six dimensional real Lie algebras contains some
errors and some omissions). We must restrict ourselves to certain special cases that
are of interest, either for matematical or physical reasons.

\bigskip 

Table 1 shows the Lie algebras of dimension $\leq8$ with non-trivial
Levi decomposition and having no invariants. Due to the low dimenions, the
only Levi factors that appear are the simple Lie algebras $\frak{so}\left(
3\right)  $ and $\frak{sl}\left(  2,\mathbb{R}\right)  $. These algebras are
of own interest, since they play an important role in multidimensional cosmologies \cite{Tu3}.%

\begin{table}
\caption{\label{t1} Lie algebras with non-trivial Levi factor and $\mathcal{N}=0$.}
\begin{tabular}{@{}|l|l|l|l|}\hline
Levi factor $\frak{s}$ & dim & {\rm Representation}
 & {\rm Nonzero structure constants}\\\hline
$\frak{sl}\left(2,\mathbb{R}\right)$ & $6$ & $D_{\frac{1}{2}}\oplus D_{0}$  & $ C_{12}^{2}=2,C_{13}^{3}=-2,C_{23}^{1}=1,C_{14}^{4}=1,C_{15}^{5}=-1$\\
 & & & $C_{25}^{4}=1,C_{34}^{5}=1,C_{46}^{4}=1,C_{56}^{5}=1$\\\hline
$\frak{so}\left(3\right)  $ & $8$ & $R_{4}\oplus D_{0}$ & $
C_{12}^{3}=1,C_{13}^{2}=-1,C_{23}^{1}=1,C_{14}^{7}=\frac{1}{2},C_{15}%
^{6}=\frac{1}{2}$\\
 & & & $ C_{16}^{5}=-\frac{1}{2},C_{17}^{4}=-\frac{1}{2},C_{24}^{5}=\frac{1}{2}%
,C_{25}^{4}=-\frac{1}{2},C_{26}^{7}=\frac{1}{2}$\\
 & & & $ C_{27}^{6}=-\frac{1}{2},C_{34}^{6}=\frac{1}{2},C_{35}^{7}=-\frac{1}{2},C_{36}^{4}=-\frac{1}{2},C_{37}^{5}=\frac{1}{2}$\\
 & & & $ C_{48}^{4}=1,C_{58}^{5}=1,C_{68}^{6}=1,C_{78}^{7}=1$\\\hline
$\frak{so}\left(  3\right)  $ & $8$ & $R_{4}\oplus D_{0}$ & $
C_{12}^{3}=1,C_{13}^{2}=-1,C_{23}^{1}=1,C_{14}^{7}=\frac{1}{2},C_{15}%
^{6}=\frac{1}{2}$\\
 & & & $C_{16}^{5}=-\frac{1}{2},C_{17}^{4}=-\frac{1}{2},C_{24}^{5}=\frac{1}{2}%
,C_{25}^{4}=-\frac{1}{2},C_{26}^{7}=\frac{1}{2}$\\
 & & & $C_{27}^{6}=-\frac{1}{2},C_{34}^{6}=\frac{1}{2},C_{35}^{7}=-\frac{1}{2},C_{36}^{4}=-\frac{1}{2},C_{37}^{5}=\frac{1}{2}$\\
 & & & $C_{48}^{4}=p,C_{48}^{6}=-1,C_{58}^{5}=p,C_{58}^{7}=-1,C_{68}^{4}=1$\\
 & & & $C_{68}^{6}=p,C_{78}^{5}=1,C_{78}^{7}=p$\\\hline
$\frak{sl}\left(  2,\mathbb{R}\right)  $ & $8$ & $2D_{\frac{1}{2}}\oplus
D_{0}$ & $C_{12}^{2}=2,C_{13}^{3}=-2,C_{23}^{1}=1,C_{14}^{4}=1,C_{15}^{5}=-1$\\
 & & & $C_{16}^{6}=1,C_{17}^{7}=-,C_{25}^{4}=1,C_{27}^{6}=1,C_{34}^{5}=1$\\
 & & &  $C_{36}^{7}=1,C_{48}^{4}=1,C_{58}^{5}=1,C_{68}^{4}=1,C_{68}^{6}=1$\\
 & & & $C_{78}^{5}=1,C_{78}^{7}=1$\\\hline
$\frak{sl}\left(  2,\mathbb{R}\right)  $ & $8$ & $2D_{\frac{1}{2}}\oplus
D_{0}$ & $C_{12}^{2}=2,C_{13}^{3}=-2,C_{23}^{1}=1,C_{14}^{4}=1,C_{15}^{5}=-1$\\
 & & & $C_{16}^{6}=1,C_{17}^{7}=-,C_{25}^{4}=1,C_{27}^{6}=1,C_{34}^{5}=1$\\
 & & & $C_{36}^{7}=1,C_{48}^{4}=1,C_{58}^{5}=1,C_{68}^{6}=p,C_{78}^{7}=p$
\\\hline
$\frak{sl}\left(  2,\mathbb{R}\right)  $ & $8$ & $2D_{\frac{1}{2}}\oplus
D_{0}$ & $C_{12}^{2}=2,C_{13}^{3}=-2,C_{23}^{1}=1,C_{14}^{4}=1,C_{15}^{5}=-1$\\
 & & & $C_{16}^{6}=1,C_{17}^{7}=-1,C_{25}^{4}=1,C_{27}^{6}=1,C_{34}^{5}=1$\\
 & & & $C_{36}^{7}=1,C_{48}^{4}=p,C_{48}^{6}=-1,C_{58}^{5}=p,C_{58}^{7}=-1$\\
 & & & $C_{68}^{4}=1,C_{68}^{6}=p,C_{78}^{5}=1,C_{78}^{7}=p$\\\hline
$\frak{sl}\left(  2,\mathbb{R}\right)  $ & $8$ & $D_{\frac{1}{2}}\oplus3D_{0}$%
& $C_{12}^{2}=2,C_{13}^{3}=-2,C_{23}^{1}=1,C_{14}^{4}=1,C_{15}^{5}=-1$\\
 & & & $C_{25}^{4}=1,C_{34}^{5}=1,C_{46}^{4}=1,C_{56}^{5}=1,C_{78}^{8}=1$\\\hline
$\frak{sl}\left(  2,\mathbb{R}\right)  $ & $8$ & $D_{\frac{3}{2}}\oplus D_{0}$
& $C_{12}^{2}=2,C_{13}^{3}=-2,C_{23}^{1}=1,C_{14}^{4}=3,C_{15}^{5}=1$\\
 & & & $C_{16}^{6}=-1,C_{17}^{7}=-3,C_{25}^{4}=3,C_{26}^{5}=2,C_{27}^{6}=1$\\
 & & & $C_{34}^{5}=1,C_{35}^{6}=2,C_{36}^{7}=3,C_{48}^{4}=1,C_{58}^{5}=1$\\
& & & $C_{68}^{6}=1,C_{78}^{7}=1$\\\hline
\end{tabular}
\end{table}

We convene that the term $D_{J}$ denotes the real representation of $\frak{sl}\left(
2,\mathbb{R}\right)  $ in its standard form, while $R_{4}$ denotes the four
dimensional real irreducible representation of $\frak{so}\left(  3\right)  $ and $D_{0}$ 
denotes the trivial representation in both cases.

\medskip

Although the general classification of these algebras seems not realizable, since it is based
on the possibility of classifying the solvable Lie algebras,
once an example is known we can deduce the following generic result:

\begin{theorem}
Let $\frak{s}$ be a semisimple Lie algebra and $\left(  R,\frak{r}\right)  $
be a pair formed by a representation of $\frak{s}$ and a solvable Lie algebra
$\frak{r}$ such that $\mathcal{N}\left(  \frak{s}\overrightarrow{\oplus}%
_{R}\frak{r}\right)  =0$. Then, for any $k\geq1$ there exists a Lie algebra
$\frak{g}_{k}$ with Levi factor $\frak{s}$ and dimension $n=\dim\left(
\frak{s}\overrightarrow{\oplus}_{R}\frak{r}\right)  +2k$ such that
$\mathcal{N}\left(  \frak{g}\right)  =0$.
\end{theorem}

\begin{proof}
Consider the Lie algebra $\frak{g}=\frak{s}\overrightarrow{\oplus}_{R^{\prime
}}\frak{r}^{\prime}$, where $R^{\prime}=R\oplus2kD_{0}$ and the radical is $\frak{r}^{\prime
}=\frak{r\oplus} k\frak{r}$ $_{2}$, where $\frak{r}_{2}$ is the affine Lie
algebra generated by $Y,Z$ and brackets $\left[  Y,Z\right]  =Z.$ The algebra
$\frak{r}^{\prime}$ is obviously a $\frak{s}$-module, and since $\mathcal{N}\left(  \frak{s}\overrightarrow{\oplus}%
_{R}\frak{r}\right)  =0$ and $\mathcal{N}\left(  \frak{r}_{2}\right)  =0$, the assertion follows from lemma 1.
\end{proof}

\begin{corollary}
For any dimension $n\geq6$ there exist Lie algebras $\frak{g}$ with
non-trivial Levi decomposition such that $\mathcal{N}\left(  \frak{g}\right)
=0$.
\end{corollary}

This reduces the classification to the pairs $\left(  R,\frak{r}\right)  $
formed by a representation of $\frak{s}$ (this being fixed) and radicals $\frak{r}$
which are indecomposable, i.e., that do not decompose into a direct sum of
ideals. Even for low dimensions like ten, it is far from being easy to find
such pairs. As an example, consider the representation $R=R_{4}\oplus3D_{0}$
of $\frak{s0}\left(  3\right)  $ and the radical $\frak{r}$ defined by the
brackets
\begin{equation*}
\begin{array}
[c]{l}%
\left[  X_{i},X_{8}\right]    =X_{i},\;4\leq i\leq7\\
\left[  X_{4},X_{9}\right]    =X_{6},\;\left[  X_{5},X_{9}\right]
=X_{7},\;\left[  X_{6},X_{9}\right]  =-X_{4}\\
\left[  X_{7},X_{9}\right]    =-X_{5},\;\left[  X_{9},X_{10}\right]  =X_{10}%
\end{array}
\end{equation*}
over the basis $\left\{  X_{4},..,X_{10}\right\}  $. This is the simplest
non-decomposable solvable Lie algebra such that the semidirect sum
$\frak{so}\left(  3\right)  \overrightarrow{\oplus}_{R}\frak{r}$ has no
non-trivial invariants (for the considered representation). In fact more is true, namely the 
nonexistence of solvable Lie algebras $\frak{r}$ such that the action of the generators $X\in\frak{r}-[\frak{r},\frak{r}]$ over the nilradical $[\frak{r},\frak{r}]$ is diagonal. This will happen also for other representations different from the one taken here. 

\section{Levi factors $\frak{s}=\frak{so}\left(  3\right)  ,\frak{sl}\left(
2,\mathbb{R}\right)  $}

Theorem 1 is a general result which holds for any Lie algebra satisfying $(8)$, and therefore non dependent on the particular Levi factor taken. Now an inspection of table 1
points out some interesting facts for the considered \ Levi factors
$\frak{s}=\frak{s0}\left(  3\right)$ and $\frak{sl}\left(  2,\mathbb{R}\right)  $.
In this section we analyze the semidirect sums $\frak{s}\overrightarrow
{\oplus}_{R}\frak{r}$ with these Levi subalgebras in more detail. Through this
section, and unless otherwise stated, the notation $\frak{s}$ will refer
either to $\frak{sl}\left(  2,\mathbb{R}\right)  $ or to $\frak{so}\left(
3\right)  $.

\medskip

We saw in section 2 that in the computation of the invariants of the algebra
$\frak{so}\left(  3\right)  \overrightarrow{\oplus}_{R}\frak{r}$ with
$R=ad\,\frak{so}\left(  3\right)  $ and $\frak{r}$ the three dimensional
abelian algebra $3L_{1}$ there was an invariant depending only on the variables
associated to $3L_{1}$. We claimed that the existence of this
 invariant, coming from a
special subsystem of $\left(  {5}\right)  $, was not casual. The next
proposition shows that this property does not depend on the representation:

\begin{theorem}
Let $R$ be an irreducible representation of $\frak{s}$. Then the semidirect
sum $\frak{s}\overrightarrow{\oplus}_{R}\frak{r}$ \ admits non-trivial
invariants. \newline Moreover, if dim$\left(  \frak{r}\right)  >\dim\left(
\frak{s}\right)  $, there exists a fundamental set of invariants formed by
functions $F_{i}$ depending only on variables associated to elements of
$\frak{r.}$
\end{theorem}

\begin{proof}
We prove it for $\frak{s}=\frak{sl}\left(  2,\mathbb{R}\right)  $,
the case of $\frak{s0}\left(  3\right)  $ being similar. At first, we only
need to prove the result for odd dimensional representations $D_{j}$, since
the remaining case follows at once from the odd dimensionality of the
semidirect sum. By proposition 1, the radical $\frak{r}$ is abelian, and the
maximal weight of $R$ is $\lambda=2m-4\;\left(  m\geq3\right)  .$ Let
$\left\{  X_{1},X_{2},X_{3},..,X_{2m}\right\}  $ be a basis of $\frak{s}%
\overrightarrow{\oplus}_{R}\frak{r}$ such that $\left\{  X_{1},X_{2}%
,X_{3}\right\}  $ is a basis of $\frak{sl}\left(  2,\mathbb{R}\right)  $ (with
$\left[  X_{1},X_{2}\right]  =2X_{2},\left[  X_{1},X_{3}\right]
=-2X_{3},\left[  X_{2},X_{3}\right]  =X_{1}$) and $\left\{  X_{4}%
,..,X_{2m}\right\}  $ a basis of the abelian radical $\frak{r}$. The system of
PDEs giving the invariants of $\frak{s}\overrightarrow{\oplus}_{R}\frak{r}$
is:%
\begin{equation}
\left.
\begin{array}
[c]{l}%
\widehat{X}_{1}F=\left(  -2x_{2}\partial_{x_{2}}+2x_{3}\partial_{x_{3}}%
-\sum_{i=0}^{2m-4}\left(  \lambda-2i\right)  x_{4+i}\partial_{x_{4+i}}\right)
F=0\\
\widehat{X}_{2}F=\left(  2x_{2}\partial_{x_{1}}-x_{1}\partial_{x_{3}}%
-\sum_{i=1}^{2m-4}\left(  \lambda-i+1\right)  x_{3+i}\partial_{x_{4+i}%
}\right)  F=0\\
\widehat{X}_{3}F=\left(  -2x_{3}\partial_{x_{1}}+x_{1}\partial_{x_{2}}%
-\sum_{i=0}^{2m-5}\left(  i+1\right)  x_{5+i}\partial_{x_{4+i}}\right)  F=0\\
\widehat{X}_{4+i}F=\left(  \left(  \lambda-2i\right)  x_{4+i}\partial_{x_{1}%
}-\left(  i+1\right)  x_{5+i}\partial_{x_{2}}+\left(  \lambda-i+1\right)
x_{3+i}\partial_{x_{3}}\right)  F=0,\;\\
0\leq i\leq2m-4
\end{array}
\right\}  ,
\end{equation}

Observe that since $\frak{r}$ is abelian, the equations $\left\{  \widehat
{X}_{4+i}F=0\;\right\}  _{0\leq i\leq2m-4}$ do not involve the partial
derivatives $\partial_{x_{i}}F$ for $4\leq i\leq2m.$ This allows us to extract
the subsystem:
\begin{equation}
\left.
\begin{array}
[c]{l}%
\widehat{X}_{1}^{\prime}F=\left(  \sum_{i=0}^{2m-4}\left(  \lambda-2i\right)
x_{4+i}\partial_{x_{4+i}}\right)  F=0\\
\widehat{X}_{2}^{\prime}F=\left(  \sum_{i=1}^{2m-4}\left(  \lambda-i+1\right)
x_{3+i}\partial_{x_{4+i}}\right)  F=0\\
\widehat{X}_{3}^{\prime}F=\left(  \sum_{i=0}^{2m-5}\left(  i+1\right)
x_{5+i}\partial_{x_{4+i}}\right)  F=0
\end{array}
\right\}  ,
\end{equation}

and any solution is obviously an invariant of $\frak{s}\overrightarrow{\oplus
}_{R}\frak{r}$. The question reduces to show that the system $\left(
{10}\right)  $ admits a non-trivial solution for any irreducible representation
$D_{J}$. Observe that $\left(  {10}\right)  $ can be written as
\begin{equation}
\left(
\begin{array}
[c]{ccccc}%
\lambda x_{4} & \left(  \lambda-2\right)  x_{5} & .. & -\left(
\lambda-2\right)  x_{2m-1} & -\lambda x_{2m}\\
0 & \lambda x_{4} & .. & 2x_{2m-2} & x_{2m-1}\\
x_{5} & 2x_{6} & .. & \lambda x_{2m} & 0
\end{array}
\right)  \left(
\begin{array}
[c]{c}%
\partial_{x_{4}}F\\
.\\
.\\
\partial_{x_{2m}}F
\end{array}
\right)  =0.
\end{equation}

Now this matrix of coefficients has at most rank three (indeed three if
$m\geq4$ and rank one if $m=3$), so that $\left(  {10}\right)  $ has always a
solution, which shows that $\mathcal{N}\left(  \frak{s}\overrightarrow{\oplus
}_{R}\frak{r}\right)  \neq0$. In particular the system $\left(  {10}\right)  $
gives the following number of solutions:
\begin{equation}
\left.
\begin{array}
[c]{l}%
1\quad {\rm if }\quad m=3\\
2m-6\quad  {\rm if}\quad m\geq4
\end{array}
\right\}  ,
\end{equation}

Observe that for $m=3$ the representation $R$ is the adjoint representation,
and in this case we can find another invariant which depends also on the
variables $x_{1},x_{2},x_{3}$. For $m\geq4$ is it not difficult to see that
$\partial_{x_{i}}F=0$ for $i=1,2,3$, which shows that the $\left(
2m-6\right)  \;$\ functionally independent solutions of $\left(  {10}\right)  $
constitute a fundamental set of invariants for $\frak{s}\overrightarrow
{\oplus}_{R}\frak{r}$.
\end{proof}

\begin{corollary}
Let $s=\frak{sl}\left(  2,\mathbb{R}\right)  ,\frak{so}\left(  3\right)  $. If
the radical $\frak{r}$ is abelian then $\mathcal{N}\left(  \frak{s}%
\overrightarrow{\oplus}_{R}\frak{r}\right)  \neq0$.
\end{corollary}

\begin{proof}
If the representation contains a copy of the trivial representation $D_{0}$ or
$\dim\frak{r}$ \ is even , we automatically have solutions of the corresponding system
$\left(  {2}\right)  $. If $R$ does not contain a copy of $D_{0}$, we can again
extract a subsystem from $\left(  {2}\right)  $, since the radical is abelian
and its equations do not contain the partial derivatives corresponding to
elements of $\frak{r}$. Now $R$ is a sum of irreducible representations, of
which at least one summand $R_{0}$ must have even highest weight $\lambda$, in
order to ensure the odd dimensionality of $\frak{r}$. Moreover, the variables
involved in $R_{0}$ do not appear in the other summands of $\frak{r}$, which
ensures that we can apply the preceding theorem. This shows that there exists
a nontrivial of the subsystem corresponding to $R_{0}$, which, by the
complete reducibility of $R$ and the abelianity of $\frak{r}$, is also an
invariant of $\frak{s}\overrightarrow{\oplus}_{R}\frak{r}$.
\end{proof}

\bigskip
 
The following example illustrates the procedure used in this proof: Let
$\frak{s}=\frak{sl}\left(  2,\mathbb{R}\right)  $ and consider the reducible
representation $D_{1}\oplus D_{\frac{1}{2}}$. Suppose that the radical
$\frak{r}$ is a five dimensional abelian Lie algebra. The invariants of
$\frak{s}\overrightarrow{\oplus}_{R}\frak{r}$ are the solutions of the system:%
\begin{equation}
\left.
\begin{array}
[c]{l}%
\left(  -2x_{2}\partial_{x_{2}}+2x_{3}\partial_{x_{3}}-2x_{4}\partial_{x_{4}%
}+2x_{6}\partial_{x_{6}}-x_{7}\partial_{x_{7}}+x_{8}\partial_{x_{8}}\right)
F=0\\
\left(  -2x_{2}\partial_{x_{1}}+x_{1}\partial_{x_{3}}+2x_{4}\partial_{x_{5}%
}+x_{5}\partial_{x_{6}}+x_{7}\partial_{x_{8}}\right)  F=0\\
\left(  2x_{3}\partial_{x_{1}}-x_{1}\partial_{x_{2}}+x_{5}\partial_{x_{4}%
}+2x_{6}\partial_{x_{5}}+x_{8}\partial_{x_{7}}\right)  F=0\\
\left(  -2x_{4}\partial_{x_{1}}-x_{5}\partial_{x_{3}}\right)  F=0\\
\left(  -2x_{4}\partial_{x_{2}}-2x_{6}\partial_{x_{3}}\right)  F=0\\
\left(  2x_{6}\partial_{x_{1}}-x_{5}\partial_{x_{2}}\right)  F=0\\
\left(  -x_{7}\partial_{x_{1}}-x_{8}\partial_{x_{3}}\right)  F=0\\
\left(  x_{8}\partial_{x_{1}}-x_{7}\partial_{x_{2}}\right)  F=0
\end{array}
\right\}  ,
\end{equation}

We extract a subsystem from the first three equations:%
\begin{equation}
\left.
\begin{array}
[c]{l}%
\left(  -2x_{4}\partial_{x_{4}}+2x_{6}\partial_{x_{6}}-x_{7}\partial_{x_{7}%
}+x_{8}\partial_{x_{8}}\right)  F=0\\
\left(  2x_{4}\partial_{x_{5}}+x_{5}\partial_{x_{6}}+x_{7}\partial_{x_{8}%
}\right)  F=0\\
\left(  x_{5}\partial_{x_{4}}+2x_{6}\partial_{x_{5}}+x_{8}\partial_{x_{7}%
}\right)  F=0
\end{array}
\right\}  ,
\end{equation}
and any solution of this system is an invariant of the algebra. $\left(
{14}\right)  $ can also be reduced to
\begin{equation}
\left.
\begin{array}
[c]{l}%
\left(  -2x_{4}\partial_{x_{4}}+2x_{6}\partial_{x_{6}}\right)  F=0\\
\left(  2x_{4}\partial_{x_{5}}+x_{5}\partial_{x_{6}}\right)  F=0\\
\left(  x_{5}\partial_{x_{4}}+2x_{6}\partial_{x_{5}}\right)  F=0
\end{array}
\right\}  ,
\end{equation}
which is the subsystem corresponding to the adjoint representation. Clearly
the polynomial $I_{1}=4x_{4}x_{6}-x_{5}^{2}$ is
a solution of $\left(  {14}\right)  $ and $\left(  {15}\right)  $, and
therefore an invariant of the algebra. Since the other summand of $R$ is
$D_{\frac{1}{2}}$, the other invariant will depend on all the variables
$x_{4},..,x_{8}$. We find $I_{2}=x_{4}x_{8}^{2}-x_{5}x_{7}x_{8}+x_{6}x_{7}%
^{2}$. Thus $I_{1},I_{2}$ form a fundamental set of invariants of $\frak{s}\overrightarrow
{\oplus}_{R}\frak{r}$.

\bigskip

These two preceding results constitute an important restriction for a
semidirect sum $\frak{s}\overrightarrow{\oplus}_{R}\frak{r}$ to satisfy
$\mathcal{N}\left(  \frak{s}\overrightarrow{\oplus}_{R}\frak{r}\right)  =0$.
Any representation in such an algebra must be reducible and contain a copy of
the trivial representation $D_{0}$ (see table 1 and the examples in section 3).

\begin{proposition}
Let $\frak{s}=\frak{so}\left(  3\right)  ,\frak{sl}\left(  2,\mathbb{R}%
\right)  $. If the radical $\frak{r}$ of $\frak{s}\overrightarrow{\oplus}%
_{R}\frak{r}$ has a one dimensional centre, then the representation $R$
describing the semidirect sum contains a copy of the trivial representation
$D_{0}$. In particular, $\mathcal{N}\left(  \frak{s}\overrightarrow{\oplus
}_{R}\frak{r}\right)  \neq0$.
\end{proposition}

\begin{proof}
\ Let $z$ generate the centre $Z\left(  \frak{r}\right)  $ or $\frak{r}$. For
any $X\in\frak{s}$ and $Y\in\frak{r}$ we have
\[
\left[  X,\left[  Y,Z\right]  \right]  +\left[  Z,\left[  X,Y\right]  \right]
+\left[  Y,\left[  Z,X\right]  \right]  =0,
\]
which shows that $\left[  \frak{s},Z(\frak{r})\right]  \subset Z\left(
r\right)  $. Now%
\[
\left[  X_{2},\left[  X_{3},Z\right]  \right]  +\left[  Z,\left[  X_{2}%
,X_{3}\right]  \right]  +\left[  X_{3},\left[  Z,X_{2}\right]  \right]  =0,
\]
which shows that $\left[  X_{1},Z\right]  =0$. Similarly it is proven that
$\left[  X_{2},Z\right]  =\left[  X_{3},Z\right]  =0$, from which we deduce
the existence of a copy of the trivial representation in the decomposition of
$R$. Since the action of $\frak{s}$ over $Z\left(  \frak{r}\right)  $ is zero,
we will obtain the monomial invariant $I_{1}=z$.
\end{proof}

The results obtained so far for the Levi factors $\frak{so}\left(  3\right)  $
and $\frak{sl}\left(  2,\mathbb{R}\right)  $ have important physical
applications, like the classification of multidimensional spacetimes \cite{De}. In this
frame, all ten dimensional real Lie algebras having a $\left(  7+d\right)
$-dimensional compact subalgebra have been determined. Of special interest are
those which have non-trivial Levi decomposition, and which are the only
candidates which could present the anomaly $\mathcal{N}\left(  \frak{s}%
\overrightarrow{\oplus}_{R}\frak{r}\right)  =0$. From the thirty classes found
\cite{Tu}, only six are indecomposable, i.e., they do not decompose as a
direct sum of lower dimensional Lie algebras. They have been listed in table
2, where the notation for the algebras is the same as in \cite{Tu3}:%

\begin{table}
\caption{\label{t2} Ten dimensional indecomposable Lie algebras with a compact subalgebra of dimension $n\geq 7$.}
\begin{tabular}{@{}|l|l|l|l|}\hline
Algebra & Levi decomposition & Representation $R$ & $\mathcal{N}$\\\hline
$L_{10,14}$ & $\frak{so}\left(  3\right)  \overrightarrow{\oplus}_{R}\left(
7L_{1}\right)  $ & $R_{7}$ & 4 \\\hline
$L_{10,15}$ & $\frak{so}\left(  3\right)  \overrightarrow{\oplus}_{R}\left(
7L_{1}\right)  $ & $R_{4}\oplus ad\,\frak{so}\left(  3\right)  $ & 4\\\hline
$L_{10,27}$ & $\frak{sl}\left(  2,\mathbb{R}\right)  \overrightarrow{\oplus
}_{R}\left(  7L_{1}\right)  $ & $D_{3}$ & 4\\\hline
$L_{10,28}$ & $\frak{sl}\left(  2,\mathbb{R}\right)  \overrightarrow{\oplus
}_{R}\left(  7L_{1}\right)  $ & $D_{2}\oplus D_{\frac{1}{2}}$ & 4\\\hline
$L_{10,29}$ & $\frak{sl}\left(  2,\mathbb{R}\right)  \overrightarrow{\oplus
}_{R}\left(  7L_{1}\right)  $ & $D_{\frac{3}{2}}\oplus D_{1}$ & 4\\\hline
$L_{10,30}$ & $\frak{sl}\left(  2,\mathbb{R}\right)  \overrightarrow{\oplus
}_{R}\left(  7L_{1}\right)  $ & $D_{1}\oplus2D_{\frac{1}{2}}$ & 4\\\hline
\end{tabular}
\end{table}

By theorem 2 and corollary 2 we see that, since the radical is always abelian, we
will obtain non-trivial invariants. For these algebras, in contrast to the possible multidimensional cosmological models seen in section 3 and table 1, the existence of a compact subalgebra of dimension $n\geq 7$ implies 
that the algebra has non-vanishing invariants. 

\section{Application to radicals with a codimension one abelian ideal}

In this section we analyze a special kind of radicals. We will suppose that
$\frak{r}$ is a solvable non-nilpotent Lie algebra such that $\left[
\frak{r},\frak{r}\right]  $ is a codimension one abelian ideal. We will see that such radicals always imply the existence of invariants, up to the lower dimensional cases. In particular the radicals found in table 1 for the 8 dimensional algebras $\frak{s}\overrightarrow{\oplus}_{R}\frak{r}$ will constitute the exception for radicals of this type. 

\begin{theorem}
Suppose that $R=R^{\prime}\oplus2D_{0}$, where $R^{\prime}$ is a
representation of $\frak{s}$. Then $\mathcal{N}\left(  \frak{s}\overrightarrow
{\oplus}_{R}\frak{r}\right)  >0$.
\end{theorem}

\begin{proof}
Since $R$ contains at least two copies of the trivial representation, there
exists an element $Y\in\left[  \frak{r},\frak{r}\right]  $ such that $\left[
\frak{s},Y\right]  =0$. Let $T\notin\left[  \frak{r},\frak{r}\right]  $ and
$\left[  T,Y\right]  =\sum_{Y_{i}\in\left[  \frak{r},\frak{r}\right]  }%
a_{i}Y_{i}\;\left(  a_{i}\in\mathbb{R}\right)  $. The equation $\widehat
{Y}F=0$ of system $\left(  {2}\right)  $ has the form%
\begin{equation}
\widehat{Y}F=-\left(  \sum_{Y_{i}\in\left[  \frak{r},\frak{r}\right]  }%
a_{i}y_{i}\right)  \partial_{T}F=0.
\end{equation}
Now, if $\left[  T,Y\right]  =0$, the function $F=y$ is an invariant of
$\frak{s}\overrightarrow{\oplus}_{R}\frak{r}$. If the bracket $\left[
T,Y\right]  $ is nonzero, then $\left(  {16}\right)  $ implies that
$\partial_{T}F=0$ for any invariant $F$. The complete reducibility of the
representation $R$ \ (the ideal $\left[  \frak{r},\frak{r}\right]  $ \ has
codimension one in $\frak{r}$ and is an $\frak{s}$-module) implies that
$\left[  \frak{s},T\right]  =0$. The number $\mathcal{N}\left(  \frak{s}%
\overrightarrow{\oplus}_{R}\frak{r}\right)  $ is given by the difference of
the dimension of $\frak{s}\overrightarrow
{\oplus}_{R}\frak{r}$ and the rank of the matrix $A\left(  \frak{s}\overrightarrow
{\oplus}_{R}\frak{r}\right)  $, which in this case has the form:%
\begin{equation}
\left(
\begin{array}
[c]{cccccccc}%
0 & \left[  X_{1},X_{2}\right]  & \left[  X_{1},X_{3}\right]  & \left[
X_{1},Z_{1}\right]  & ... & \left[  X_{1},Z_{r}\right]  & 0 & 0\\
\left[  X_{2},X_{1}\right]  & 0 & \left[  X_{2},X_{3}\right]  & \left[
X_{2},Z_{1}\right]  & ... & \left[  X_{2},Z_{r}\right]  & 0 & 0\\
\left[  X_{3},X_{1}\right]  & \left[  X_{3},X_{2}\right]  & 0 & \left[
X_{3},Z_{1}\right]  & ... & \left[  X_{3},Z_{r}\right]  & 0 & 0\\
\left[  Z_{1},X_{1}\right]  & \left[  Z_{1},X_{2}\right]  & \left[
Z_{1},X_{3}\right]  & 0 & ... & 0 & 0 & \left[  T,Z_{1}\right] \\
\vdots & \vdots & \vdots & \vdots &  & \vdots & \vdots & \vdots\\
\left[  Z_{r},X_{1}\right]  & \left[  Z_{r},X_{2}\right]  & \left[
Z_{r},X_{3}\right]  & 0 & ... & 0 & 0 & \left[  T,Z_{r}\right] \\
0 & 0 & 0 & 0 & ... & 0 & 0 & \left[  T,Y\right] \\
0 & 0 & 0 & \left[  Z_{1},T\right]  & ... & \left[ Z_{r},T\right]  &
\left[  Y,T\right]  & 0
\end{array}
\right)  ,
\end{equation}

where $\left\{  Z_{1},..,Z_{r},Y,T\right\}  $ is a basis of $\frak{r}$.
Elementary methods show that the determinant of this matrix is the product of $-\left[  T,Y\right]^{2}$ and the following determinant:
\begin{equation}
\det\left(
\begin{array}
[c]{cccccc}%
0 & \left[  X_{1},X_{2}\right]  & \left[  X_{1},X_{3}\right]  & \left[
X_{1},Z_{1}\right]  & ... & \left[  X_{1},Z_{r}\right] \\
\left[  X_{2},X_{1}\right]  & 0 & \left[  X_{2},X_{3}\right]  & \left[
X_{2},Z_{1}\right]  & ... & \left[  X_{2},Z_{r}\right] \\
\left[  X_{3},X_{1}\right]  & \left[  X_{3},X_{2}\right]  & 0 & \left[
X_{3},Z_{1}\right]  & ... & \left[  X_{3},Z_{r}\right] \\
\left[  Z_{1},X_{1}\right]  & \left[  Z_{1},X_{2}\right]  & \left[
Z_{1},X_{3}\right]  & 0 & ... & 0\\
\vdots & \vdots & \vdots & \vdots &  & \vdots\\
\left[  Z_{r},X_{1}\right]  & \left[  Z_{r},X_{2}\right]  & \left[
Z_{r},X_{3}\right]  & 0 & ... & 0
\end{array}
\right),
\end{equation}
which must be zero, since the rank of the matrix in $\left(  {18}\right)  $
gives the number of invariants of the subalgebra $\frak{s}\overrightarrow
{\oplus}_{R-2D_{0}}\left(  rL_{1}\right)  $, which is non-maximal in virtue of
theorem 1. Therefore the rank of $A\left(  \frak{s}\overrightarrow{\oplus}%
_{R}\frak{r}\right)  $ is less than its dimension, from which the existence of
non-trivial invariants is ensured.
\end{proof}

\medskip

It should be remarked that if $R$ contains only one copy of $D_{0}$ or the
codimension of $\left[  \frak{r},\frak{r}\right]  $ is $\frak{r}$ is greater
than one, then the conclusion is false, as can easily be extracted from table
1. We will finally see that radicals as considered in this section are only
valid in low dimensions in order to obtain Lie algebras $\frak{s}%
\overrightarrow{\oplus}_{R}\frak{r}$ such that $\mathcal{N}\left(
\frak{s}\overrightarrow{\oplus}_{R}\frak{r}\right)  =0$.

\begin{proposition}
If $\dim\left(  \frak{r}\right)  \geq7$ then $\mathcal{N}\left(
\frak{s}\overrightarrow{\oplus}_{R}\frak{r}\right)  \neq0$.
\end{proposition}

\begin{proof}
Like before, since $\left[  \frak{r},\frak{r}\right]  $ is a codimension one
$\frak{s}$-submodule of $\frak{r}$, the action of $\frak{s}$ on a generator
$T\in\frak{r}-\left[  \frak{r},\frak{r}\right]  $ is zero. If $\dim\left(
\frak{r}\right)  =7$, then $\dim\left(  \frak{s}\overrightarrow{\oplus}%
_{R}\frak{r}\right)  =10$ and the matrix $A\left(  \frak{s}\overrightarrow
{\oplus}_{R}\frak{r}\right)  $ has the form%
\begin{equation}
\left(
\begin{array}
[c]{ccccccc}%
0 & \left[  X_{1},X_{2}\right]  & \left[  X_{1},X_{3}\right]  & \left[
X_{1},Z_{1}\right]  & ... & \left[  X_{1},Z_{6}\right]  & 0\\
\left[  X_{2},X_{1}\right]  & 0 & \left[  X_{2},X_{3}\right]  & \left[
X_{2},Z_{1}\right]  & ... & \left[  X_{2},Z_{6}\right]  & 0\\
\left[  X_{3},X_{1}\right]  & \left[  X_{3},X_{2}\right]  & 0 & \left[
X_{3},Z_{1}\right]  & ... & \left[  X_{3},Z_{6}\right]  & 0\\
\left[  Z_{1},X_{1}\right]  & \left[  Z_{1},X_{2}\right]  & \left[
Z_{1},X_{3}\right]  & 0 & ... & 0 & \left[  T,Z_{1}\right] \\
\vdots & \vdots & \vdots & \vdots &  & \vdots & \vdots\\
\left[  Z_{6},X_{1}\right]  & \left[  Z_{6},X_{2}\right]  & \left[
Z_{6},X_{3}\right]  & 0 & ... & 0 & \left[  T,Z_{6}\right] \\
0 & 0 & 0 & \left[  Z_{1},T\right]  & ... & \left[  Z_{6},T\right]  & 0
\end{array}
\right)  .
\end{equation}

It is routine to verify that the determinant of $\left(  {19}\right)  $ does not
depend on the brackets, and that it is zero. Since for any radical $\frak{r}$
of the considered type such that $\dim\left(  \frak{r}\right)  \geq7$ the
determinant of $A\left(  \frak{s}\overrightarrow{\oplus}_{R}\frak{r}\right)  $
$\ $is a linear combination of matrices of type $\left(  {19}\right)  $ and
matrices like in $\left(  {18}\right)  $, it follows that $\det A\left(
\frak{s}\overrightarrow{\oplus}_{R}\frak{r}\right)  =0$.
\end{proof}

\medskip
Observe that this result explains, in terms of the representation theory of $\frak{so}(3)$, why the ten dimensional Galilei algebra has two (Casimir) invariants depending only on the translations $P_{i}$ and the pure Galilean transformations $K_{i}$.  

\section{Conclusions}

We have seen that for any dimension $n\geq 6$ there exist non-semisimple Lie
algebras $\frak{s}\overrightarrow{\oplus }_{R}\frak{r}$ with non-trivial
Levi factor $\frak{s}$ and such that $\mathcal{N}\left( \frak{s}%
\overrightarrow{\oplus }_{R}\frak{r}\right) =0$. This constitutes a proof
that the Levi decomposition theorem \cite{RWZ} does not reduce the number of
generalized Casimir invariants of $\frak{s}\overrightarrow{\oplus }_{R}\frak{%
r}$ to some combination of the numbers corresponding to the Levi factor $%
\frak{s}$ and the radical $\frak{r}$, but depends essentially on the pair $%
\left( R,\frak{r}\right) $ \ formed by the representation $R$ describing the
semidirect sum and the radical. 

\medskip

For the rank one simple Lie algebras $\frak{so}\left( 3\right) $ and $\frak{%
sl}\left( 2,\mathbb{R}\right) $ the number of invariants of a semidirect sum 
$\frak{s}\overrightarrow{\oplus }_{R}\frak{r}$ have been analized in some
detail. In particular, the analysis undertaken has given a representation
theoretic interpretation of the invariants obtained for the $\left(
3+1\right) $ kinematical algebras like the Galilei algebra. The interest of
these Levi factors is therefore justified not only by kinematical problems,
but also by the extensions of Bianchi type-IX cosmology \cite{Tu,De}.
Specially interesting are those admissible extensions which have no
invariants. Therefore invariant quantities for these algebras should be searched using distribution theory \cite{PWZ}. 
In particular, if the radical is abelian, we have proved that we
will obtain solutions, some of them depending only on variables associated
to elements of the radical. This confirms that the fact that the special
affine algebras $\frak{sa}\left( n,\mathbb{R}\right) $ have invariants (for
being odd dimensional) is not an isolated case, but also the general pattern
for those semidirect sums which are even dimensional. From the computed examples, it seems reasonable to expect that, whenever the radical $\frak{r}$ is a nilpotent Lie algebra, the number of invariants of a semidirect sum $\frak{s}\overrightarrow{\oplus }_{R}\frak{r}$ will be nonzero. However, for this case it is not sufficient to know which is the representation $R$ that describes the semidirect sum. We need more precise information on the structure of $\frak{n}$ (not merely the value of very general invariants like the nilpotence index), which impedes to establish a general result as for the abelian case.

The most important question that arises from our results is whether they can
be extended to any semisimple Lie algebra of rank $r\geq 2$. At least for
direct sums of $\frak{sl}\left( 2,\mathbb{R}\right) $ and $\frak{so}\left(
3\right) $ this seems to hold. An example which is worth to be analyzed is
the Schr\"{o}dinger algebra $\frak{S}$ in $\left( 3+1\right) $ dimensions 
\cite{Nie}. Over the basis $\left\{ J_{i},K_{i},P_{i},P_{0},C,D\right\}
_{i=1,2,3}$ this algebra is given by the brackets:

\begin{equation*}
\begin{tabular}{llll}
$\left[ J_{i},J_{j}\right] =\varepsilon _{ijk}J_{k},$ & $\left[ J_{i},K_{j}%
\right] =\varepsilon _{ijk}K_{k},$ & $\left[ J_{i},P_{j}\right] =\varepsilon
_{ijk}P_{k},$ & $\left[ K_{i},P_{0}\right] =P_{i},$ \\ 
$\left[ P_{i},D\right] =P_{i},$ & $\left[ D,K_{j}\right] =D_{j},$ & $\left[
D,P_{0}\right] =-2P_{0},$ & $\left[ C,P_{j}\right] =K_{j},$ \\ 
$\left[ C,P_{0}\right] =-D,$ & $\left[ C,D\right] =-2C.$ &  & 
\end{tabular}
\end{equation*}

\smallskip

where $P_{0}$ is the time translation, $P_{i}$ the space translations, $J_{i}
$ the rotations and $K_{i}$ the pure Galilean transformations. It can easily
be verified that the subalgebra $\frak{a}$ \ generated by $\left\{
K_{i},P_{i}\right\} _{i=1,2,3}$ is six dimensional and abelian, while $%
\left\{ P_{0},C,D\right\} $ generates a copy pf $\frak{sl}\left( 2,\mathbb{R}%
\right) $. Therefore we obtain the semisimple algebra $\frak{so}\left(
3\right) \oplus \frak{sl}\left( 2,\mathbb{R}\right) $, and since $\frak{a}$
is an ideal, we have the Levi decomposition of $\frak{S}$ (by abuse of notation we can denote the corresponding representation by $D_{\frac{1}{2}}\otimes ad \frak{so}(3)$). If we extract a
subsystem of the corresponding  system $\left( {2}\right) $, as done in the proof of theorem 2, we obtain that $\frak{S}$ has a fourth order Casimir operator $%
\mathcal{P}_{4}$ depending only on the space translations and pure Galilei
transformations:
\begin{eqnarray*}
\mathcal{P}_{4} &=&K_{1}^{2}\left( P_{2}^{2}+P_{3}^{2}\right)
+K_{2}^{2}\left( P_{1}^{2}+P_{3}^{2}\right) +K_{3}^{2}\left(
P_{1}^{2}+P_{2}^{2}\right)  \\
&&-2\left(
P_{1}P_{2}K_{1}K_{2}+P_{1}P_{3}K_{1}K_{3}+P_{2}P_{3}K_{2}K_{3}\right)  
\end{eqnarray*}

\bigskip For other simple Lie algebras a direct calculation of the rank of
matrices $A\left( \frak{s}\overrightarrow{\oplus }_{R}\frak{r}\right) $
becomes a enormously difficult problem, and therefore the proofs of the
generalization of the results obtained for rank one simple algebras, if they hold, must be
approached by completely different means. 

Finally, these results are of interest for the study of non-semisimple
(maximal) regular subalgebras of simple Lie algebras. The example $\frak{sl}%
\left( 2,\mathbb{R}\right) \overrightarrow{\oplus }_{D_{\frac{1}{2}}\oplus D_{0}}3L_{1}%
$ of table 1 is a regular subalgebra of $\frak{sl}\left( 3,%
\mathbb{R}\right) $ and has no invariants. It would be important to obtain a detailed description of the non-semisimple maximal regular subalgebras of simple Lie algebras which do not have non-trivial invariants. This problem is of interest not only for symmetry
breaking questions \cite{Dy}, but also for solving many fundamental problems
which arise in rigidity and contraction theory \cite{Cam,We}, like the
invariant theory of parabolic subalgebras of semisimple Lie algebras, the
construction of contraction trees or the expansion problem \cite
{AL}.

\end{document}